\documentclass{amsart}
\usepackage{graphicx}
\theoremstyle{definition}

\theoremstyle{remark}

\numberwithin{equation}{section}
\begin{document}

\title{Some Nice Sums are Almost as Nice if you turn them Upside Down}
\maketitle

\vskip 10pt

\begin{center}
{\bf Moa Apagodu}\\
{Department of Mathematics, Virginia Commonwealth University,\\ Richmond, VA 08854, USA}

\vskip 10pt

{\bf Doron Zeilberger} \\
Department of Mathematics, Rutgers University,
Hill Center-Busch Campus, \\110 Frelinghuysen Rd., Piscataway,
NJ 08854-8019, USA.
\end{center}

\begin{abstract}

  {We represent the sums $\sum_{k=0}^{n-1}{n \choose k}^{-2}\,\,,\sum_{k=0}^{m}{m\choose k}^{-1}{a\choose n-k}^{-1},$  $\sum_{k=0}^{n-1}\frac{q^{-k(k-1)}}{{\genfrac{[}{]}{0pt}{}{n}{k}}_q }$,
and the sum of the reciprocals of the summands in Dixon's identity, each as a product of an {\it indefinite hypergeometric sum}
times a (closed form) {\it hypergeometric sequence}} .

\end{abstract}

\vskip 10pt

\noindent In [3], Rockett proved the formula
$$
\sum_{k=0}^{n}{n \choose k}^{-1}  = \frac {n+1}{2^n}\sum_{j=0}^{n}\frac{2^j}{j+1}\,\,, \eqno(Rockett)
$$

\noindent for all nonnegative integers $n$. Since then,
several authors (e.g. [1][4][5]) used different techniques to re-derive (\emph{Rockett}) and considered
analogous, more complicated, sums.\\

As we all know,
$$
\sum_{k=0}^{n}{n \choose k}  = 2^n \,\,, \eqno(Binomial)
$$
is the {\it simplest} example of a {\bf definite} binomial coefficient sum that evaluates in {\it closed-form}.
Such sums have the form
$$
\sum_{k=0}^{n} F(n,k)  = A(n) \,\, ,
$$
where $A(n)$ is {\it hypergeometric} (i.e. $A(n+1)/A(n)$ is a rational function of $n$),
and $F(n,k)$ is {\it bi-hypergeometric} (i.e $F(n+1,k)/F(n,k)$  and $F(n,k+1)/F(n,k)$ are both rational functions of $n$ and $k$).
Other famous examples, are
$$
\sum_{k=0}^{n} {n \choose k}{a \choose k}  = {n+a \choose a} \,\, ,
\eqno(Chu-Vandermonde)
$$
and
$$
\sum_{k=-n}^{n} (-1)^k {a+n \choose n+k}  {b+n \choose b+k}  {a+b \choose a+k}
=
\frac{(a+b+n)!}{a!b!n!}
 \,\, .
\eqno(Dixon)
$$
The {\it best} thing that can happen to a binomial coefficients sum is to be evaluable in {\it closed form},
i.e. as a hypergeometric sequence. The {\it next best} thing is to be ``almost closed-form'', i.e. of the form
$$
ClosedForm (n) \cdot \left ( \sum_{j=1}^{n} h(j) \right ) \quad ,
$$
where $h(j)$ is hypergeometric in the {\it single} variable $j$.

Inspired by $(Rockett)$, that has this ``almost perfect'' form, and that came out from considering
the sum of the reciprocals of the summands of the {\it simplest} binomial coefficient sum known,
we searched for other cases where one starts from a well-known binomial coefficients sum that
evaluates in closed form, and looks at what happens if one considers the sum of the reciprocals
of the summands. To our great surprise and delight, it worked for the Chu-Vandermonde summand and
for the Dixon summand (see below). To our disappointment, it didn't work for the so-called
Pfaff-Saalschutz identity (see, e.g. [2]).

\noindent  We use the WZ method([2][7][8]), and the reader is assumed to be familiar with it.
In particular we used Zeilberger's Maple packages {\tt EKHAD} (procedure zeillim) and {\tt qEKHAD} accompanying
[2], available from

{\tt http://www.math.rutgers.edu/\~{}zeilberg/tokhniot/[EKHAD,qEKHAD]} .

\bigskip

\noindent In most applications of the WZ method, the summation limits are {\it natural},
i.e. the summand is well-defined and vanishes at $k=-1$ and $k=n+1$ (or in the case of $(Dixon)$,
at $k=-n-1$ and $k=n+1$)
so one can get a {\it homogeneous} first-order linear recurrence for the definite sum,
that leads to a closed-form solution. In the present cases,
this is no longer the case, and we only get
{\it inhomogeneous} first-order linear recurrences, that lead to the
above-mentioned ``almost-perfect'' kind of solutions.

\noindent Consider the sum
$$
f(n):=\sum_{k=\alpha}^{n-\beta}F(n,k)\,\,,
$$

\noindent where $\alpha$ and $\beta$ are constants and $[\alpha,n-\beta]$ is properly contained in $[0,n]$, the support of $F(n,k)$.
Assume also that there exists a function $G(n,k)$, (the so-called certificate) such that
$$
a(n)F(n+1,k)+b(n)F(n,k)=G(n,k+1)-G(n,k)\,\,. \eqno(WZeqn)
$$

\noindent Take $\alpha=0$ and $\beta=1$. Then, if we add both sides of (\emph{WZeqn}) from $k=0$ to $k=n$ and rewrite the result we get a nonhomogeneous recurrence relation satisfied by $f(n)$:

$$
a(n)f(n+1)+b(n)f(n)= G(n,n)-G(n,0)+a(n)F(n+1,n) \,\,  . \eqno{(nonhomorec)}
$$

\noindent Finally, we solve (\emph{nonhomorec}) and get the following representation of $f(n)$ as an indefinite sum:
$$
f(n)=g(n)\sum_{i=0}^{n-1}\frac{C(i)}{g(i)a(i)}\,\,,
$$
\noindent
where
$$
C(i)= G(i,i)-G(i,0)+a(i)F(i+1,i) \,\,  .
$$
(i.e. the right side of $(nonhomorec)$ with $n$ replaced by $i$),
and $g(n)$ is a solution to the associated homogeneous equation $a(n)f(n+1)+b(n)f(n)=0$. \\

\noindent We note that in case $\alpha \neq 0$ and $\beta\neq 1$, the left-hand side of (\emph{nonhomrec}) remains the same. The only change is on the right-hand side. \\

\noindent \textbf{Theorem 1 [Upside Down Binomials Squared]}:

$$
\sum_{k=0}^{n-1}{n \choose k}^{-2}=\frac{(n+1)!(n+1)^2}{(n+\frac{3}{2})!4^n}\sum_{j=0}^{n-1}\frac{(3j^3+12j^2+18j+10)4^j(j+3/2)!}
{(j+1)^2(j+1)!(j+2)^3} \,\,.
$$

\noindent (Here $n!:=\Gamma(n+1)$).

\noindent \textbf{Proof}: We construct the function

$$
G(n,k)=\frac{(2k-3n-6)(k-n-1)^2}{{n \choose k}^{2}}
$$

\noindent such that $(4n+10)(n+1)^2F(n+1,k)-(n+2)^2F(n,k)=G(n,k+1)-G(n,k)$, where $F(n,k)$ is the summand on the left-hand side sum.\\

\noindent Next add both sides from $k=0$ to $k=n$ and rearrange to get the nonhomogeneous recurrence relation satisfied by the sum on the left-hand side, $S_n$ :

$$
(4n+10)(n+1)^2S_{n+1}-(n+2)^2S_n=3n^3+12n^2+18n+10\,\,\,.
$$

\noindent Finally the theorem follows by solving this recurrence with the initial condition $S_1=1$.\\

\noindent Next we give q-analog of (\emph{Rockett}). Let $(a;q)_n= (1-a)(1-aq)\ldots(1-aq^{n-1})$. \\

\vskip .2in

\noindent \textbf{Theorem  2 [Upside Down q-binomial]} :

$$
\sum_{k=0}^{n-1} \frac{q^{-k(k-1)/2}}{{\genfrac{[}{]}{0pt}{}{n}{k}}_q  } = \frac{(q^2;q)_n}{(q^2;q^2)_n}\sum_{i=0}^{n-1}
\frac{C(i)(q^2;q^2)_i}{(q^2;q)_i(q^{i+2}-1)}\,\,,
$$

\noindent where
$$C(i)=\frac{q^{i+1}+q^{3(i+1)}+q^{-i(i-1)/2}+q^{-(i+1)(i-4)/2}-q^{(2+3i-i^2)/2}-q^{-(i+1)(i-2)/2}-2q^{2(i+1)}}{q^{i+1}-1}\,\,.$$

\noindent \textbf{Proof}: We construct the function

$$
G(n,k)=\frac{q^{n+1}(1-q^{n-k+1})q^{-k(k-1)/2}}{ {\genfrac{[}{]}{0pt}{}{n}{k}}_q }
$$

\noindent such that $(q^{n+1}-1)(q^{n+1}+1)F(n+1,k)-(q^{n+2}-1)F(n,k)= G(n,k+1)-G(n,k)$, where $F(n,k)$ is the summand on the left-hand side.\\

\noindent Add both sides from $k=0$ to $k=n$ and rearrange to get the nonhomogeneous recurrence relation satisfied by the sum on the left-hand side, $S_n$ :\\

$$
(q^{n+1}-1)(q^{n+1}+1)S_{n+1}-(q^{n+2}-1)S_n = C(n)\,\,\,.
$$

\noindent Finally the theorem follows by solving this recurrence with the initial condition $S_1=1$. \\

\vskip .2in

\noindent Next we consider the reciprocal of the summand in the Chu-Vandermonde and the Dixon classical identities. \\

\vskip .2in

\noindent \textbf{Theorem 3 [Upside Down Chu-Vandermonde]}  \\

$$
\sum_{k=0}^{m-1}{m\choose k}^{-1}{a\choose n-k}^{-1}=g(m)\sum_{i=0}^{m-1}\frac{C(i)}{g(i)(i+2)(a+i-n+2)}
$$

\noindent where

$$
g(m)= \frac{(a+m-n+1)!(a+4)!(m+1)}{2(a+m+3)!(a-n+2)!}
$$

and\\

$$
C(i)=(n+i+in+1){a\choose n}^{-1} + (2i+a-n+3){a\choose n-i}^{-1}\,\,.
$$

\noindent \textbf{Proof}: The proof follows similarly from the recurrence relation

$$
(a+m+4)(m+1)S_{m+1}-(m+2)(a+m-n+2)S_m = C(m)\,\,.
$$

\noindent where $S_m$ is the sum on left-hand side.

\vskip .2in

\noindent \textbf{Theorem 4 [Upside Down Dixon]}  \\

$$
\sum_{k=0}^{n-1}(-1)^k{n+b\choose n+k}^{-1}{n+c\choose c+k}^{-1}{b+c\choose b+k}^{-1} = g(n)\sum_{i=0}^{n-1}\frac{C(i)}{g(i)(2i+2)(c+i+2)(b+i+2)}
$$

\noindent where\\

$$
g(n)= \frac{(b+n+1)(c+n+1)(b+c+2)!n!}{(b+1)(c+1)(b+c+n+2)!}
$$

\noindent and\\

$
C(i)=(-1)^i(3bi+3ci+bc+3b+5i^2+12i+3c+7){i+b\choose 2i}^{-1}{b+c\choose b+i}^{-1} -
(b+1)(c+1)(i+1)
{i+b\choose i}^{-1}{i+c\choose c}^{-1}{b+c\choose b}^{-1}\,\,.
$\\

\vskip .2in

\noindent \textbf{Proof}: The proof follows similarly from the recurrence relation

$$
2(c+n+1)(b+n+1)(b+c+n+3)S_{n+1}-2(n+1)(c+n+2)(b+n+2)S_n = C(n)\,\,.
$$

\noindent where $S_n$ is the sum on left-hand side.

\vskip .2in

\noindent \textbf{Remarks} a) From the recurrence in the proof of theorem 1,

$$
2S_n=\frac{(n+1)^3}{2n^3+3n^2}S_{n-1}+\frac{3n^3+3n^2+3n+1}{2n^3+3n^2}
$$

\noindent that implies

$$
\lim_{n\rightarrow \infty} \sum_{k=0}^{n}{n \choose k}^{-2}=2\,\,.
$$

\noindent b) Identity (\emph{Rockett}) is a special case of the following identity with $x=y=1$ in

$$
\sum_{k=0}^{n}{n \choose k}^{-1}x^ky^{n-k}=x^n+\left(\frac{xy}{x+y}\right)^n(n+1)\sum_{j=0}^{n-1}\frac{((j+1)y^{j+2}+yx^{j+1})(x+y)^j}
{(xy)^{j+1}(j+1)(j+2)}
$$

\noindent which follows from the recurrence [$S_n$ the left-hand side sum]:\\

$$
(x+y)(n+1)S_{n+1}-(n+2)xyS_n = (n+1)y^{n+2}+yx^{n+1}\,\,.
$$

\vskip .2in

\noindent c) T. Mansour [1] derived the following representation for the reciprocals of binomial coefficients to any power \emph{m} as:

$$
\sum_{k=0}^{n-1}{n \choose k}^{-m}=(n+1)^m\sum_{k=0}^{n}\left[\sum_{j=0}^{k}\frac{(-1)^j}{{n-k+1+j}} {k \choose j}\right]^m\,\,.
$$

\noindent One of the referees kindly pointed out that (b) and (c) can also be deduced from a general identity derived in [6].
We thank them for this remark and for numerous corrections and improvements.
\\

\noindent {\bf References}\\

\bibliographystyle{amsplain}

\noindent [1]  T. Mansour. Combinatorial Identities and Inverse Binomial Coefficients. \emph{Advances in Applied Mathematics}, \textbf{28}: 196-202, 2002.\\

\noindent [2]  M. Petkovsek, H.S. Wilf, D. Zeilberger. A=B. A.K. Peters Ltd., 1996.\\

\noindent [3] A.M. Rockett. Sum of the inverse of binomial coefficients.  \emph{Fibonacci Quart.}, 19: 433-437, 1981.\\

\noindent [4] R. Sprugnoli. Sums of Reciprocals of the Central Binomial Coefficients. \emph{Integers} \textbf{6}: \#A27,2006.\\

\noindent [5] B. Sury. Sum of the Reciprocals of the Binomial Coefficients. \emph{Europ. J. Combinatorics}, 14: 351-353, 1993. \\

\noindent [6] B. Sury,  T. Wang, and F. Zhao. Identities involving reciprocals of binomial coefficients. \emph{Journal of Integer Sequences}, Vol.7: Article 04.2.8, 2007. \\

\noindent [7]  H. Wilf, D. Zeilberger. Rational Functions Certify Combinatorial Identities. \emph{J. Amer. Math. Soc}. 3: 147-158, 1990. \\

\noindent [8]  D. Zeilberger. The method of creative telescoping. \emph{J. Symbolic Computation}, 11: 195-204, 1991.

\end{document}